\documentclass[12pt, a4paper]{amsart}
\usepackage[T1]{fontenc}
\usepackage{amsfonts}
\usepackage{listings}
\usepackage{color}
\usepackage{fancyvrb}
\usepackage{amsmath}
\usepackage{mathrsfs}
\usepackage{amsthm}
\usepackage{graphicx}
\usepackage{enumerate}
\usepackage{bbm}
\usepackage{euscript}
\usepackage{yhmath}
\usepackage{MnSymbol}
\usepackage{parskip}
\usepackage[colorlinks=true,
            linkcolor=blue,  
            citecolor=blue,  
            urlcolor=blue]{hyperref}
\usepackage[a4paper, margin=1in]{geometry}
\newcommand{\R}{\mathbb{R}}

\newcommand{\ep}{\varepsilon}

\newcommand{\mbf}{\mathbf}

\newcommand{\LL}{\mathrm{L}}
\newcommand{\WW}{\mathrm{W}}

\newcommand{\dd}{\mathrm{d}}
\newcommand{\BMO}{\mathrm{BMO}}

\newtheorem{theorem}{Theorem}[section]
\newtheorem{lemma}[theorem]{Lemma}
\newtheorem{prop}[theorem]{Proposition}

\theoremstyle{definition}
\newtheorem{definition}[theorem]{Definition}

\DeclareMathOperator{\dist}{dist}

\DeclareMathOperator{\Div}{div}

\title[Asymptotic regularity for variational integrals]
{Sharp regularity results for variational
integrals asymptotic to the $p$-Dirichlet integral}
\author{Aidan Strong}
\begin{document}
\baselineskip=18pt

\maketitle

\begin{abstract}
We consider the regularity for minimisers of $\int_\Omega F(Du) \dd x$ where $F(z) \approx |z|^p$ in suitable senses at infinity. Depending on the strength of this asymptotic condition, we discuss various regularity results, and by means of counterexamples we show they are sharp. An effort is made to obtain results without assuming unnecessary regularity of $F$ while covering the full range $1<p<\infty$.
\end{abstract}
\setcounter{tocdepth}{1}
\begingroup
\small
\tableofcontents
\endgroup
\section{Introduction}
We consider local minimisers of the variational integral 
$$\mathcal{F}[u]=\int_\Omega F(Du) \dd x$$
where $\Omega \subset \R^n$ is a bounded open set, and $F:\R^{N \times n} \rightarrow \R$ is a continuous function henceforth called an integrand. We consider the following question: Given  $1<p<\infty$ and an integrand $F(z)$ which is asymptotic to $|z|^p$ in a suitable sense, what is the regularity of minimisers of $\int_\Omega F(Dv) \dd x$? This is a classical question, firstly sparked from the famous result of Uhlenbeck \cite{Uhlenbeck1977} who showed that $p$-harmonic maps, which correspond to $F(z)=|z|^p$, are locally $C^{1,\alpha}$ (see Proposition \ref{prop p harmonic functions}). The first result in the asymptotic setting was due to Chipot--Evans \cite{ChipotEvans1986}, who proved that if  $F\in C^2(\R^{N \times n})$ and $$\lim_{|z| \rightarrow \infty}F''(z) = A$$ for some inner product $A$ on $\R^{N \times n}$, then minimisers of $\mathcal{F}$ are locally Lipschitz. The case $A=2I$ here represents asymptotic behaviour to the Dirichlet integral in a $C^2$ sense. Later that year, Giaquinta--Modica \cite{GiaquintaModica1986}
 extended this result to the case $p \geq 2$ by proving minimisers are locally Lipschitz whenever $F \in C^2(\R^{N \times n})$ satisfies
$$\lim_{|z| \rightarrow \infty}\frac{|F''(z)-D^2(|z|^p)|}{|z|^{p-2}}=0.$$

If instead one considers $F$ asymptotic to $|z|^p$ in a zeroth-order sense, namely that $$\lim_{|z| \rightarrow \infty} \frac{|F(z)-|z|^p|}{|z|^p}=0,$$
it was shown by Fuchs \cite{Fuchs1994} that $u \in C_\mathrm{loc}^{0,\alpha}$ for all $\alpha<1$, and later Dolzmann-Kristensen \cite{DolzmannKristensen2005} showed that $u \in \WW^{1,q}_\mathrm{loc}$ for all $q<\infty$. Similar results were also shown to hold for extremals, namely solutions of the Euler--Lagrange equation, in Kristensen--Taheri \cite{KristensenTaheri2003}.
We also note the work of Scheven--Schmidt \cite{SchevenSchmidt2010}, who prove higher integrability results in a more general case where $F$ is asymptotic to some suitably strongly convex function $G$. In all of the above results, no quasiconvexity of the integrand $F$ is assumed. The standard methods of partial regularity, where one compares $u$ to a solution of a linearised problem at small scales, are therefore not available. Instead, the arguments typically involve comparing $u$ to suitable $p$-harmonic maps. \\ 
In this paper we attempt to unify the presentation and treatment of the above results. Alongside streamlining the presentation, there are four main contributions in this work. Firstly, all of our results cover the full range $1<p<\infty$, which is achieved by heavily relying on the auxiliary maps $V_p(z)=|z|^{\frac{p-2}{2}}z$ and $A_p(z)=|z|^{p-2}z$. Secondly, we consider the case where $F(z) \approx |z|^p$ in a suitable first-order sense, alongside the zeroth-order and second-order asymptotic conditions that have previously been studied. Thirdly, we relax the regularity conditions on $F$. Finally, by means of counterexamples we show that all of our results are sharp. Most of our proofs are based on direct comparison with $p$-harmonic maps, as opposed to blow-up arguments. It is worth noting that the integrands considered here give a large class of functionals in the vectorial setting $N \geq  2$ where everywhere regularity holds.  The main results presented in this paper are the following:
\begin{enumerate}[(i)]
    \item If $F \in C(\R^{N \times n})$ is a bounded perturbation of $|z|^p$ then $V_p(Du) \in \BMO_\mathrm{loc}$ (Theorem \ref{theorem BMO in the case of extremely close}). However, $u$ is not necessarily Lipschitz (Theorem \ref{theorem counterexample non lipschitz extremely close}).
    \item If $F \in C(\R^{N \times n})$ is asymptotic to $|z|^p$  in a zeroth-order sense then $u \in \WW^{1,q}_\mathrm{loc}$ for all $q<\infty$ (Theorem \ref{thm higher integrability close to dirichlet}). However, $Du$ is not necessarily $\BMO$ (Theorem \ref{theorem counterexample dirichlet at infinity C^0})
    \item If $F \in C(\R^{N \times n}) $ is asymptotic to $|z|^p$ in a first-order sense then  $Du \in \BMO_\mathrm{loc}$ (Theorem \ref{Theorem BMO regularity C^1 approximation for all p}), however $u$ is not necessarily Lipschitz.
    \item If $F \in C^1(\R^{N \times n})$ is asymptotic to $|z|^p$ in a second-order sense then $u \in \WW^{1,\infty}_\mathrm{loc}$ (Theorem \ref{theorem lipschitz regularity under approximate C^2 sense}), however $Du$ is not necessarily continuous.
\end{enumerate} 
Asymptotic in a zeroth-order sense simply means $|F(z)-|z|^p| = o(|z|^p)$. In the case $F \in C^1(\R^{N \times n})$, our definition of asymptotic in a first-order sense is 
$$\limsup_{|z| \rightarrow \infty} \frac{|F'(z)-pA_p(z)|}{|z|^{p-2}}<\infty.$$
In the case $F \in C^2(\R^{N \times n})$, our definition of asymptotic in a second-order sense agrees with the one considered in \cite{GiaquintaModica1986}, namely that $$\lim_{|z| \rightarrow \infty}\frac{|F''(z)-D^2(|z|^p)|}{|z|^{p-2}}=0.$$ 
We of course recall that Theorem \ref{thm higher integrability close to dirichlet} is a special case of results from \cite{DolzmannKristensen2005}, however we include a short proof in our setting for completeness. The proof of Theorem \ref{theorem BMO in the case of extremely close} is in some sense a simpler version of Theorem \ref{thm higher integrability close to dirichlet}, although we do not know a precise reference. Our proof of Theorem \ref{Theorem BMO regularity C^1 approximation for all p} relies on $\BMO$ estimates for degenerate equations in \cite{ {Diening2012}}. Finally, the proof of Theorem \ref{theorem lipschitz regularity under approximate C^2 sense} is similar to classical methods, but we carefully adapt the proof to work for all $1<p<\infty$ without assuming $F \in C^2$. For our counterexamples, we construct integrands $F \in C(\R^{2 \times 2})$ using the conformal splitting of matrices, inspired by the papers \cite{DolzmannKristensenZhang2013, dolzmann2012}. Our original motivation for revisiting this problem came from studying the regularity of minimisers for the distance functional $\int_{\Omega}\mathrm{dist}(Du,K)^p \dd x$ where $K \subset \R^{N \times n}$ is a compact set, which is useful in applications to multi-well problems. The integrand $F_K^p(z)=\dist(z,K)^p$ need not be $C^2$, but under suitable geometric conditions on $K$, it is asymptotic to $|z|^p$ in the second order sense as described in this paper. More generally, the relaxation of this functional is given by $\int_\Omega (F_K^p)^{qc}(Du) \dd x$, where $qc$ denotes the quasiconvex envelope. It turns out that $(F_K^p)^{qc}$ is asymptotic to $|z|^p$ in the first-order sense as described in this paper. Details of these facts are addressed in our companion paper \cite{Strong2026}.

\section{Preliminaries}
\subsection{Notation}
Throughout this paper, $\Omega \subset \R^n$ denotes an open bounded set where $n \geq 2$. We consider maps $u \in \WW^{1,p}(\Omega;\R^{N})$ where $1<p<\infty$ and write $\R^{N \times n}$ for the space of $N\times n$ matrices. We equip $\mathbb R^{N\times n}$ with the standard
Euclidean inner product denoted $z \cdot w$, and denote the associated norm by $|z|$. We write
$p'=p/(p-1)$ for the conjugate exponent of $p$. For $x_0\in\mathbb R^n$ and $R>0$, we denote by
$B_R(x_0)$ the open ball of radius $R$ centred at
$x_0$, and as is standard we simply write $B_R$ when it is convenient. For an integrable function $f$, we denote the mean over the ball $B_R(x_0)$ by $$(f)_R=(f)_{B_R(x_0)}=\strokedint_{B_R}f \dd x = \frac{1}{|B_R|}\int_{B_R}f \dd x.$$
As is standard, $C,c>0$ may refer to constants changing from line to line. We call a constant $C$ universal if it only depends on the data $n,N,p$. We write $A\lesssim B$ if $A\leq CB$ for some constant $C>0$, and $A \approx  B$ if $A\lesssim B$ and $B\lesssim A$. We sometimes also use big $O$ and little $o$ notation. Given an integrand $F:\R^{N\times n}\rightarrow  \R$ with $p$-growth, we say that $u\in \WW^{1,p}_{\mathrm{loc}}
(\Omega;\mathbb R^N)$ is a local minimiser of
$\int_\Omega F(Dv)\dd x$ if
$$\int_{\Omega'} F(Du) \dd x
\leq \int_{\Omega'} F(Du+D\varphi)\dd x$$
for every $\Omega' \Subset \Omega$ and
$\varphi\in W^{1,p}_0(\Omega';\R^N)$.
\subsection{$p$-harmonic maps}
\begin{definition}
Let $1<p<\infty$ and $\Omega \subset \R^n$ be an open set. We say $h \in \WW^{1,p}(\Omega;\R^N)$ is $p$-harmonic if any of the following equivalent definitions hold:
\begin{enumerate}[(i)]
    \item $-\Div (|Dh|^{p-2}Dh)=0$ in the sense of distributions on $\Omega$
    \item $\int_{\Omega} |Dh|^{p-2}Dh \cdot D\varphi \dd x =0$ for all $\varphi \in W_0^{1,p}(\Omega;\R^N)$
    \item $h$ is a local minimiser of $\int_{\Omega}|Du|^p \dd x$.
\end{enumerate}
\end{definition}
It is convenient to work with the auxiliary functions
$$A_p(z)=|z|^{p-2}z,$$
$$V_p(z)=|z|^{\frac{p-2}{2}}z,$$
where $z,w \in \R^{N \times n}$.
Note that $|z|^p$ is a $C^1$ function whose derivative is precisely $pA_p(z).$ Of course we understand that $A_p(0)=V_p(0)=0.$ In the following we write $(Dh)_r = \strokedint_{B_r}Dh \dd x$ and $V_p(Dh)_r = (V_p(Dh))_r$.
\begin{prop}\label{prop p harmonic functions}
Let $h \in \WW^{1,p}(\Omega;\R^N)$ be $p$-harmonic. Then there exist universal constants $\alpha \in(0,1)$ and $C>0$ such that the following hold for any concentric $B_r \subset B_R \Subset \Omega$:
\begin{enumerate}[(i)]
    \item $h \in C^{1,\alpha}_\mathrm{loc}(\Omega;\R^N)$
    \item $\strokedint_{B_r}|Dh-(Dh)_r|^2 \dd x \leq C \left(\frac{r}{R} \right)^{2\alpha} \strokedint_{B_R}|Dh-(Dh)_R|^2 \dd x $
    \item $\strokedint_{B_r}|V_p(Dh)-V_p(Dh)_r|^2 \dd x \leq C \left(\frac{r}{R} \right)^{2\alpha} \strokedint_{B_R}|V_p(Dh)-V_p(Dh)_R|^2 \dd x $
    \item $\sup_{B_{R/2}}|Dh|^p \leq C \strokedint_{B_R}|Dh|^p \dd x$.
\end{enumerate}
\end{prop}
\begin{proof}
$(i)$ is folklore originating from \cite{Uhlenbeck1977}. $(ii)$ and $(iv)$ are found in \cite[(2.4), (2.3)]{DiBenedettoManfredi1993}, and $(iii)$ is from \cite[Proposition 3]{DuzaarMingione2004}.
\end{proof}
Since many proofs rely heavily on the auxiliary functions $A_p$ and $V_p$, it is very useful to record the following equivalences. 
\begin{lemma}\label{lemma equivalences of auxiliary V functinos}
Let $1<p<\infty$ and $z,w \in \R^{N \times n}$. The following quantities are comparable up to universal constants:
\begin{enumerate}[(i)]
    \item $|V_p(z)-V_p(w)|^2$
    \item $(|z|+|w|)^{p-2}|z-w|^2$
    \item $(A_p(z)-A_p(w)) \cdot (z-w)$
    \item $|A_p(z)-A_p(w)| |z-w|$
    \item $ |z|^p -|w|^p -pA_p(w) \cdot (z-w)$
    \item $|z-w|^2 \int_0^1|tz+(1-t)w|^{p-2} \dd t$.
\end{enumerate}
\end{lemma}
\begin{proof}
We note that $(i) \approx (ii)$ is found in \cite[Lemma 1]{DuzaarMingione2004}. Moreover, the equivalences $(ii) \approx (iii) \approx (iv)$ follow from  \cite[Lemma 2.4]{diening2009}. \\
\underline{$(v) \approx (ii)$:} By the Fundamental Theorem of Calculus,
\begin{align*}
|z|^p-|w|^p-pA_p(w) \cdot (z-w) &= p \int_0^1 (A_p(w+t(z-w))-A_p(w))\cdot (z-w) \dd t \\
&\overset{(iii) \approx (ii)}{\approx }p\int_0^1 (|w+t(z-w)|+|w|)^{p-2}t \dd t |z-w|^2.
\end{align*}
So, it suffices to justify that $$\int_0^1 (|w+t(z-w)|+|w|)^{p-2} t \dd t \approx (|z|+|w|)^{p-2}.$$ This is a simple variant of a standard lemma, see \cite[Lemma 8.3]{Giusti}. Indeed, there it follows that 
$$\int_0^1t^r (\mu+|w+t(z-w)|)^s \dd t \approx (\mu+|w|+|z|)^s$$ for all $r \geq 0, s>-1$ and $\mu \geq 0$, where the implied constants are independent of $\mu$ but may depend on $s$ and $r$. Setting $\mu = |w|$, $r=1$, and $s=p-2$ gives our desired inequality.\\
\underline{$(vi) \approx (ii):$} Simply using the same integral equivalence as above, but now with $r=0$, $\mu=0$, and $s=p-2$.
\end{proof}
We record that for $p \geq 2$ the above lemma gives $|z-w|^p \leq (|z|+|w|)^{p-2}|z-w|^2 \approx|V_p(z)-V_p(w)|^2$, so that $V_p^{-1}$ is globally Hölder continuous with exponent $2/p$. 
\subsection{Dichotomy estimates}
We use the following two lemmas, which allow us to deduce higher integrability or boundedness of a map $H \in \LL^2(\Omega)$ if at each scale we have suitable control of either $\strokedint_B |H|^2 \dd x$ or the decay of the excess $\strokedint_{B}|H-(H)_B|^2 \dd x.$
\begin{lemma}\cite[Lemma 1]{DolzmannKristensen2005}\label{lemma Iwaniec higher integrability}.
Let $s>1$. There exists some constant $\sigma=\sigma_{s,n}$ such that the following holds. Let $H \in \LL^2(\Omega)$, and suppose that there exist some $\tau \in (0,1)$ and $M>0$ such that either $$\strokedint_{B} |H|^2 \dd x \leq M \qquad \text{or} \qquad \strokedint_{\tau B}|H-(H)_{\tau B}|^2 \dd x \leq \sigma \strokedint_{B}|H|^2 \dd x $$ holds for every ball $B \Subset \Omega$. Then $H \in \LL^{2s}_\mathrm{loc}(\Omega)$ and 
$$\left(\strokedint_{B_{R/2}}|H|^{2s} \dd x \right)^{1/s} \leq C\left(M+\strokedint_{B_R}|H|^2 \dd x \right)$$
for all $B_R \Subset \Omega$, where $C=C(n,\tau,s)$.
\end{lemma}
\begin{lemma}\label{lemma dichotomy boundedness}
Let $H \in \LL^2(\Omega)$, and suppose that there exist some $\tau \in(0,1)$ and $M>0$ such that either
$$\strokedint_{B}|H|^2 \dd x \leq M \qquad \text{or} \qquad \strokedint_{\tau B}|H-(H)_{\tau B}|^2 \dd x \leq \frac{1}{2}\strokedint_{B}|H-(H)_B|^2 \dd x$$
holds for every ball $B \Subset \Omega$. Then $H \in \LL^\infty_\mathrm{loc}(\Omega)$ and $$\sup_{B_{R/2}}|H|^2 \leq C \left(M+ \strokedint_{B_{R}}|H|^2 \dd x\right)$$
for every $B_R \Subset \Omega$, where $C=C(n,\tau)$.
\end{lemma}
\begin{proof}
 Fix $B_R=B_R(x_0) \Subset \Omega$ and take  $x \in B_{R/2}$ a Lebesgue point of $H$. We now consider concentric balls $B_{r_k}= B_{r_k}(x)$ where $r_k = \frac{R}{2}\tau^k$. If $\strokedint_{B_{r_k}}|H|^2 \dd x \leq M$ for infinitely many $k$, by passing to a subsequence and using the Lebesgue Differentiation Theorem, we get $|H(x)|^2 \leq M$. So, we just need to consider the other case, so that there is some (minimally chosen) $k_0$ such that for all $k\geq k_0$ we have $$E_{r_{k+1}} \leq \frac{1}{2}E_{r_k}$$
where $E_{r_k}=\strokedint_{B_{r_k}}|H-(H)_{B_{r_k}} |^2 \dd x$. Iterating this obviously gives $E_{r_k} \leq 2^{-(k-k_0)}E_{r_{k_0}}$ for all $k \geq k_0.$ Using the fact that $$|(H)_{B_{r_{k+1}}}-(H)_{B_{r_k}}| \leq \tau^{-n/2}E_{r_k}^{1/2}$$
and the Lebesgue differentiation theorem, we see that 
\begin{align*}
|H(x)|&=\lim_{k \rightarrow \infty}|(H)_{B_{r_k}}| \\
& \leq |(H)_{B_{r_{k_0}}}|+\sum_{k=k_0}^\infty|(H)_{B_{r_{k+1}}}-(H)_{B_{r_k}}|  \\
& \leq  |(H)_{B_{r_{k_0}}}|+C E_{r_{k_0}}^{1/2}\\
& \leq C \left( \strokedint_{B_{r_{k_0}}}|H|^2 \dd x\right)^{1/2}.
\end{align*}
If $k_0=0$ this gives $|H(x)|^2 \leq C \strokedint_{B_R(x_0)}|H|^2 \dd x$. Otherwise, $k_0 \neq 0$ so that $\strokedint_{B_{R_{k_0-1}}}|H|^2 \dd x \leq M$ which gives $|H(x)|^2 \lesssim M.$
\end{proof}
\subsection{Conformal splitting}
We will construct several counterexamples using the following splitting of $2 \times 2$ matrices into their conformal and anticonformal parts.
Given $X = \begin{pmatrix} X_{11} & X_{12} \\ X_{21} & X_{22} \end{pmatrix} \in \mathbb{R}^{2\times 2}$, we define 
\[ X^+ = \frac{1}{2} \begin{pmatrix} X_{11} + X_{22} & -(X_{21} - X_{12}) \\ X_{21} - X_{12} & X_{11} + X_{22} \end{pmatrix}, \quad X^- = \frac{1}{2} \begin{pmatrix} X_{11} - X_{22} & X_{21} + X_{12} \\ X_{21} + X_{12} & -(X_{11} - X_{22}) \end{pmatrix}. \]

Then $X = X^+ + X^-$ and $X^+ \cdot X^- = 0$. Moreover, $$2\det X = |X^+|^2 - |X^-|^2$$ and $$|X^+|^2 + |X^-|^2 = |X|^2.$$

\section{Bounded perturbations of the $p$-Dirichlet integral}
In this section, we consider continuous integrands $F$ which satisfy $$|F(z)-|z|^p| \leq C$$ for all $z \in \R^{N \times n}$. Equivalently, we can write  $F(z)=|z|^p+G(z)$ where $G$ is a bounded continuous function. We present this situation first, as it is the simplest example to demonstrate the direct comparison method we employ several times. The following result is possibly known to experts, but to the author's best knowledge there is no precise reference.
\begin{theorem}\label{theorem BMO in the case of extremely close}
Suppose $u \in \WW^{1,p}(\Omega;\R^N)$ is a local minimiser of $\int_{\Omega}F(Dv) \dd x$ where $F$ is a continuous integrand satisfying $|F(z)-|z|^p| \leq C$ for all $z \in \R^{N \times n}$. Then $$V_p(Du) \in \BMO_{\mathrm{loc}}(\Omega;\R^{N\times n}).$$
\end{theorem}

\begin{proof}
Fix $B_R \Subset \Omega$, and let $h \in u+\WW_0^{1,p}(B_R;\R^{N})$ be $p$-harmonic, so that $$\strokedint_{B_R}|Dh|^p \dd x \leq \strokedint_{B_R}|Du|^p \dd x.$$
On the other hand since $u$ is a minimiser, $$\strokedint_{B_R}|Du|^p \dd x \leq C+\strokedint_{B_R}F(Du) \dd x \leq C+\strokedint_{B_R}F(Dh) \dd x \leq C+\strokedint_{B_R}|Dh|^p \dd x.$$Combining these facts we have $$0 \leq \strokedint_{B_R}|Du|^p-|Dh|^p \dd x \leq C.$$ Using this, the Euler--Lagrange equation for $h$, and Lemma \ref{lemma equivalences of auxiliary V functinos}, we have $$\strokedint_{B_R}|V_p(Du)-V_p(Dh)|^2 \dd x \approx \strokedint_{B_R}|Du|^p-|Dh|^p-p|Dh|^{p-2}Dh \cdot (Du-Dh) \dd x \leq C.$$ Combining this estimate with the decay estimate for $Dh$, we can obtain the suitable decay estimate for $Du$. Indeed, let us write $$E_R = \strokedint_{B_R}|V_p(Du)-V_p(Du)_R|^2 \dd x$$ for the excess of $Du$. Then for $0<\tau <1$, using Proposition \ref{prop p harmonic functions}
\begin{align*}
E_{\tau R} & \leq \strokedint_{B_{\tau R}}|V_p(Du)-V_p(Dh)_{\tau R}|^2 \dd x \\
& \lesssim \strokedint_{B_{\tau R}}|V_p(Dh)-V_p(Dh)_{\tau R}|^2 \dd x +\strokedint_{B_{\tau R}}|V_p(Du)-V_p(Dh)|^2 \dd x \\
& \leq C \tau^{2\alpha }\strokedint_{B_{R}}|V_p(Dh)-V_p(Dh)_R|^2 \dd x + C\tau^{-n} \strokedint_{B_{ R}}|V_p(Du)-V_p(Dh)|^2 \dd x\\
& \leq C \tau^{2\alpha }\strokedint_{B_{R}}|V_p(Dh)-V_p(Du)_R|^2 \dd x + C\tau^{-n} \strokedint_{B_{ R}}|V_p(Du)-V_p(Dh)|^2 \dd x\\
& \leq C \tau^{2\alpha }E_R + C\tau^{-n}.
\end{align*}
Choosing $\theta < 2\alpha$ and $\tau$ small enough that $C\tau^{2\alpha}<\tau^\theta$, we have
 $$E_{\tau R} \leq \tau^\theta E_R+C.$$ Iterating this estimate in the standard way gives $$E_r \leq C\left(\frac{r}{R} \right)^\theta E_R+C$$ for all $r <R$. This inequality implies that $E_r$ is bounded for sufficiently small $r$. This gives $V_p(Du) \in \BMO_{\mathrm{loc}}(\Omega;\R^{N\times n})$ as claimed.
\end{proof}
Since $|V_p(Du)|^2 =|Du|^p$, it follows from the John--Nirenberg inequality that if $V_p(Du) \in \BMO_\mathrm{loc}$ then $Du \in \LL^q_\mathrm{loc}$ for any $q<\infty$. If $p \geq 2$ then $V_p^{-1}$ is globally Hölder continuous, and consequently $Du \in \BMO_\mathrm{loc}$ whenever $V_p(Du) \in \BMO_\mathrm{loc}$. By means of a counterexample, we now would like to show that Theorem \ref{theorem BMO in the case of extremely close} is sharp in the sense that we cannot obtain Lipschitz regularity without further assumptions. The following construction is based on a variant of \cite[Proposition 6.1]{DolzmannKristensenZhang2013} and \cite{dolzmann2012}. Recall that $F: \R^{2 \times 2} \rightarrow \R$ is polyconvex if $F(z)=G(z, \det z)$ for some convex function $G:\R^{(2\times 2)+1} \rightarrow \R$.

\begin{theorem}\label{theorem counterexample non lipschitz extremely close}
There exists a polyconvex integrand $F:\R^{2\times 2}\rightarrow \R$ which satisfies $$|F(z)-|z|^2| \leq C$$ such that $\int_{B}F(Dv) \dd x $ has a non Lipschitz minimiser for some ball $B$ centered at the origin.
\end{theorem}
\begin{proof}
We write $(x,y) \in \R^2$ and $r=\sqrt{x^2+y^2}$. Working on a small ball centred at the origin, we define the map $$u(x,y) = \left(\log \frac{1}{r}\right)(x,y).$$ Writing $e_r = \frac{1}{r}(x,y)$, one computes
$$Du = \left(\log\frac{1}{r}\right)I -e_r \otimes e_r$$
$$Du^+=\left(\log\frac{1}{r}-\frac{1}{2}\right)I$$
$$Du^-=\frac{1}{2}I-e_r\otimes e_r.$$
It is not difficult to see that $|Du^-|^2 = \frac{1}{2}$.\footnote{A neat way to see this: observe that $e_r$ and $e_\theta$ are eigenvectors of $Du^-$ with eigenvalues $\pm \frac{1}{2}$. The norm of a symmetric matrix is given by the Euclidean norm of vector of its eigenvalues.} Now we define the integrand $$F(X)=|X|^2 + (1-2|X^-|^2)_+.$$
Obviously $|F(X)-|X|^2| \leq 1$. Moreover writing $|X|^2 =2 \det X+2|X^-|^2$ and using the identity $t_+ = t+(-t)_+$, we can write
$$F(X)=1+2\det X + (2|X^-|^2-1)_+ \geq 1+2 \det X.$$
Now suppose that $v \in u+W_0^{1,2}(B)$. Since $\det$ is a null-Lagrangian we have $$\int_B\det Dv \dd x = \int_B \det Du \dd x.$$ On the other hand, since $|Du^-|^2=\frac{1}{2}$ we have $F(Du) = 1+2 \det Du $. Putting these together we clearly have $$\int_B F(Du) \dd x \leq \int_B F(Dv) \dd x.$$So, $u$ is a local minimiser, but clearly $Du$ is not bounded near the origin. Finally $F$ is polyconvex noting that $F(X)=\max \left\{|X|^2, 1+2\det X \right\}$.
\end{proof}
We emphasise that the key to this counterexample is constructing an integrand $F$ which takes the form $$F=\text{null Lagrangian} + \text{non-negative error}$$ and where the error part is zero when evaluated at $Du$. Our choice of $F$ above was motivated by the referenced literature, but in fact there are many different choices. Using ideas from \cite{sverakyan}, we believe it is possible to modify the above counterexample so that $F$ can be taken to be smooth and convex, however we leave the details for future work. Instead, we give another example to show that even exponential decay of $|F(z)-|z|^2|$ does not lead to Lipschitz regularity. 
\begin{theorem}
There exists a smooth integrand $F: \R^{2 \times 2}\rightarrow \R$ which satisfies $$|F(z)-|z|^2| \leq Ce^{-|z|}$$
such that $\int_B F(Dv) \dd x$ has a non Lipschitz minimiser.
\end{theorem}

\begin{proof}
We use a construction similar to the above, but this time take $u= \left(\log \log \frac{1}{r}\right)(x,y)$. We compute
$$Du=\log \log \frac{1}{r}I-\frac{1}{\log \frac{1}{r}}e_r \otimes e_r$$
$$Du^+=\left(\log \log \frac{1}{r}-\frac{1}{2 \log \frac{1}{r}} \right)I$$
$$Du^-=\frac{1}{\log \frac{1}{r}}\left(\frac{1}{2}I-e_r \otimes e_r \right).$$
One checks that 
$$|Du^-|^2=\frac{1}{2 (\log \frac{1}{r})^2}$$
and 
$$|Du|^2=2 \left(\log \log \frac{1}{r} \right)^2-\frac{2 \log \log \frac{1}{r}}{\log \frac{1}{r}}+\frac{1}{\left(\log \frac{1}{r} \right)^2}.$$
In particular, 
$|Du| \sim \sqrt 2 \log \log \frac{1}{r}$ for small $r$. 
Let $\chi: [0,\infty)\rightarrow [0,1]$ be a smooth cutoff satisfying $\chi=1$ on $[0,1]$ and $\chi=0$ on $[2,\infty)$. We define $$F(X)=|X|^2-2|X^-|^2 \chi(
|X^-|^2 e^{|X|}
)$$
which is a smooth function. One checks that $$|F(X)-|X|^2| \leq 4e^{-|X|}$$
and 
$$F(X)=2 \det X+2|X^-|^2(1-\chi(
|X^-|^2 e^{|X|}) \geq 2 \det X.$$
Working on a sufficiently small ball we have $F(Du)=2 \det Du$ since $|Du^-|^2 e^{|Du|}=\frac{1}{2}(\log \frac{1}{r})^{\sqrt 2-2}(1+o(1)) \rightarrow 0$ as $r \rightarrow 0$. Consequently $u$ is a local minimiser. 
\end{proof}
 On the other hand, it is shown in \cite{Fuchs1996} that if $F(z)=|z|^2+G(z)$ where $G$ is Lipschitz and vanishes for large $|z|$, then minimisers are in fact locally Lipschitz.

\section{Zeroth-order asymptotic regularity}
In this section, we consider integrands satisfying  $$|F(z)-|z|^p| = o(|z|^p)$$
as $|z| \rightarrow \infty$. In other words, $F(z)$ resembles $|z|^p$ at infinity in a zeroth-order sense. The arguments here are  similar to the previous section, but the difference is that we must instead settle for higher integrability rather than $\BMO$ regularity. The regularity result here is a special case of results from \cite{DolzmannKristensen2005}, however we have decided to include a short proof based on Lemma \ref{lemma Iwaniec higher integrability} for completeness and the reader's convenience. 
\begin{lemma}
Suppose $u \in \WW^{1,p}(\Omega;\R^N)$ is a local minimiser of $\int_\Omega F(Dv) \dd x $ where $F$ is a continuous integrand satisfying $|F(z)-|z|^p| = o(|z|^p)$. Let $h \in u+\WW_0^{1,p}(B_R;\R^{N})$ be $p$-harmonic. Then for every $\ep >0$ there exists $A_\ep >0$ such that $$\strokedint_{B_R}|V_p(Du)-V_p(Dh)|^2 \dd x \leq \ep \strokedint_{B_R}(|Du|^p+A_\ep) \dd x .$$
\end{lemma}

\begin{proof}
Note that our assumption means for every $\ep>0$ there exists $A_\ep >0$ such that $$(1-\ep)|z|^p -A_\ep \leq F(z) \leq (1+\ep)|z|^p + A_\ep.$$  Applying the growth bounds and minimisation, just like we did in Theorem \ref{theorem BMO in the case of extremely close}, we have
$$\strokedint_{B_R}(1-\ep)|Du|^p \dd x \leq A_\ep + (1-\ep)\strokedint_{B_R}|Dh|^p \dd x + 2\ep \strokedint_{B_R}|Du|^p \dd x.$$
This gives, up to redefining $\ep$,  $$0 \leq \strokedint_{B_R}|Du|^p-|Dh|^p \dd x \leq A_\ep +\ep\strokedint_{B_R}|Du|^p \dd x.$$ This in turn gives the claimed estimate once more following the proof of Theorem \ref{theorem BMO in the case of extremely close}.
\end{proof}

\begin{theorem}[\cite{DolzmannKristensen2005}]\label{thm higher integrability close to dirichlet}
Suppose $u \in \WW^{1,p}(\Omega;\R^N)$ is a local minimiser of $\int_\Omega F(Dv) \dd x $ where $F$ is a continuous integrand satisfying $|F(z)-|z|^p| = o(|z|^p)$. Then $u \in \WW^{1,q}_{\mathrm{loc}}(\Omega;\R^N)$ for all $q<\infty$. 
\end{theorem}
\begin{proof}
Fixing the ball $B \Subset \Omega$, it follows from the previous estimates, and the decay estimates for $p$-harmonic functions, that $$\strokedint_{\tau B}|V_p(Du)-V_p(Du)_{\tau B}|^2 \dd x \leq (C\tau^{2\alpha}+C\tau^{-n}\ep)\strokedint_{B}|V_p(Du)|^2 \dd x + C\tau^{-n}\ep A_\ep$$ for every $\ep >0$, $\tau \in (0,1)$, where $C, \alpha, A_\ep$ do not depend on the ball. Now we fix $s>1$ and let $\sigma$ denote the constant from Lemma \ref{lemma Iwaniec higher integrability}. We check the condition of Lemma \ref{lemma Iwaniec higher integrability} holds for $H = V_p(Du)$. Set $\tau \in (0,1)$ small enough that $C \tau^{2\alpha} < \sigma/3$, and then choose $\ep >0$ small enough that $C \ep \tau^{-n} < \sigma/3$. Then $$\strokedint_{\tau B}|V_p(Du)-V_p(Du)_{\tau B}|^2 \dd x \leq \frac{2\sigma}{3} \strokedint_{B}|V_p(Du)|^2 \dd x + \frac{\sigma}{3}A$$ for some constant $A$ which is independent of the ball. Now we set $M=A$. If $\strokedint_B |V_p(Du)|^2 \dd x >M=A$ we have  $\strokedint_{\tau B}|V_p(Du)-V_p(Du)_{\tau B}|^2 \dd x \leq \sigma \strokedint_{B}|V_p(Du)|^2 \dd x$ as required. By Lemma \ref{lemma Iwaniec higher integrability} we conclude that $V_p(Du) \in L^{2s}_\mathrm{loc}(\Omega;\R^{N\times n}).$ Since $s$ is arbitrary, we are done.
\end{proof} 
Note that from Sobolev embedding, the higher integrability of $Du$ implies $u \in C^{0,\alpha}_{\mathrm{loc}}$ for all $\alpha <1$. As we did in the previous section, we now want to say that Theorem \ref{thm higher integrability close to dirichlet} is sharp, this time in the sense that we cannot obtain $\BMO$ regularity of the gradient without further assumptions. The counterexample runs along the same lines as Theorem \ref{theorem counterexample non lipschitz extremely close}, however the construction of the integrand is slightly more delicate, although it still takes the form of $$F=\text{null Lagrangian} + \text{non-negative error}$$ and where the error part is zero when evaluated at $Du$, and crucially $|Du^+| \approx |Du|$ near the origin.

\begin{theorem}\label{theorem counterexample dirichlet at infinity C^0}
There exists a continuous integrand $F: \R^{2 \times 2}\rightarrow \R$ which satisfies $$|F(z)-|z|^2| = o(|z|^2)$$ as $z \rightarrow \infty$, such that $\int_{B}F(Dv) \dd x$ has a minimiser whose gradient is not $\BMO_\mathrm{loc}$, where $B=B_\rho(0)$ is a suitably small ball. 
\end{theorem}
\begin{proof}
We consider the function $u:\R^2 \rightarrow \R^2$, $$u(x,y) = |\log\frac{1}{r}|^\alpha (x,y)$$ where $r=\sqrt{x^2+y^2}$ and $\alpha>1$ is a fixed number. We work on a domain $\Omega=B_\rho$ where $\rho <1 $ is sufficiently small in a way which only depends on $\alpha$. One can compute:
\begin{align*}
Du^+ &= \left(\log\frac{1}{r}\right)^{\alpha-1} \left(\log\frac{1}{r} - \frac{\alpha}{2}\right) \begin{pmatrix} 1 & 0 \\ 0 & 1 \end{pmatrix} \\
Du^- &= - \frac{\alpha}{2r^2} \left(\log\frac{1}{r}\right)^{\alpha-1} \begin{pmatrix} x^2-y^2 & 2xy \\ 2xy & -(x^2-y^2) \end{pmatrix}.
\end{align*}
In particular, we note that (for $r$ sufficiently small)
\begin{align*}
|Du^+|&=\sqrt 2\left( \log \frac{1}{r} \right)^{\alpha-1}\left( \log \frac{1}{r}-\frac{\alpha}{2}\right)\\
|Du^-|&= \frac{\sqrt 2}{2}\alpha \left( \log \frac{1}{r}\right)^{\alpha-1} \\
|Du|&=\sqrt 2\left( \log \frac{1}{r} \right)^{\alpha} \left(1-\frac{\alpha
}{\log \frac{1}{r}}+\frac{\alpha^2}{2 \left( \log \frac{1}{r}\right)^2}\right)^{1/2}.
\end{align*}
For small $r$ we see that for some constant $c>0$, $$|\left\{x \in B_\rho: |Du(x)|>t \right\}| \geq \left|\left\{x \in B_\rho: \left(\log \frac{1}{r} \right)^\alpha > ct\right\} \right| \approx e^{- 2(ct)^{1/\alpha}}.$$ From Lemma \ref{lemma function not in BMO} this means that $Du \notin \BMO_\mathrm{loc}$.\\
\underline{Claim:} There exists a smooth function $\psi: [0,\infty) \rightarrow [0,\infty)$ such that
$$\psi(|Du^+|)=|Du^-|$$
and 
$$\psi(t) = o(t) \text{ as } t \rightarrow \infty.$$\\
\underline{Proof:} Considering our above formulae, we need that
$$\psi\left(\sqrt 2\left( \log \frac{1}{r} \right)^{\alpha-1}\left( \log \frac{1}{r}-\frac{\alpha}{2}\right) \right)=\frac{\sqrt 2}{2}\alpha \left( \log \frac{1}{r}\right)^{\alpha-1}$$
for all $r$ sufficiently small. This is obviously equivalent to 
$$\psi\left(\sqrt 2s^{\alpha-1} \left(s-\frac{\alpha}{2}\right)\right)=\frac{\sqrt 2}{2}\alpha s^{\alpha-1}$$
for all $s$ sufficiently big. But, it is clear that $s \mapsto \sqrt 2 s^{\alpha-1}\left( s-\frac{\alpha}{2}\right)$ is injective and strictly increasing for large $s$, so by taking the inverse of that map, we can define $\psi$ for large $t$, and smoothly extending $\psi$ gives our construction. Note that our construction gives $\psi(t) \approx \frac{\alpha}{2^{1-\frac{1}{2\alpha}}}t^{1-\frac{1}{\alpha}}=o(t)$ for large $t$.\\
We now complete the construction. Define $$F(X)=2 \det X+ 2(\psi(|X^+|)-|X^-|)^2$$ where $\psi$ is defined as above. Note that $$F(X)=|X|^2 - 4 \psi(|X^+|)|X^-|+2\psi(|X^+|)^2.$$ 
Since $\psi(t)=o(t)$ it follows $\psi(|X^+|) =o (|X|)$, and hence $F(X) = |X|^2 +o(|X|^2).$
Finally, $F(Du)=2 \det Du$, while if $v \in u+\WW_0^{1,2}(B_\rho,\R^N)$ then $$\int_{B_\rho}F(Dv) \dd x \geq \int_{B_\rho }2 \det Dv \dd x = \int_{B_\rho}2 \det Du \dd x$$
so $u$ is a local minimiser.
\end{proof}
Once more, we would expect it is possible to construct a counterexample where $F$ is additionally smooth and strongly convex, however we leave this for future work.

\section{First-order asymptotic regularity}
In view of the previous section, it might be natural to hope that a first order condition like $$|F'(z)-pA_p(z)| = o(|z|^{p-1})$$
would lead to some better regularity. Note that integrating this condition implies an approximation in a $C^0$ sense, namely $|F(z)-|z|^p| = o(|z|^p)$. It turns out that this condition alone is not enough to obtain $\BMO$ regularity. The following construction also gives another proof of Theorem \ref{theorem counterexample dirichlet at infinity C^0}.

\begin{theorem}\label{theroem counterexample to show no naive C^1}
There exists a smooth integrand $F:\R^{2 \times 2} \rightarrow \R$ which satisfies $$|F'(z)-2z|=o(|z|)$$
as $|z| \rightarrow \infty$, such that $\int_BF(Dv) \dd x$ has a minimiser whose gradient is not $\BMO_\mathrm{loc}$, where $B$ is a suitably small ball. 
\end{theorem}
\begin{proof}
Following the construction in Theorem \ref{theorem counterexample dirichlet at infinity C^0}, we take $u=\left|\log \frac{1}{r}\right|^\alpha(x,y)$ with $\alpha>1$ as before, and take $\psi:[0,\infty)\rightarrow [1,\infty)$ a smooth function as before satisfying $\psi(|Du^+|)=|Du^-|$, and $\psi(t)\lesssim (1+t)^{1-\frac{1}{\alpha}}$, $\psi'(t) \lesssim (1+t)^{-\frac{1}{\alpha}}$, as well as $\psi(t)=1$ near $0$. This time, we set
$$F(X)=2\det X+2 \frac{(|X^-|^2-\psi(|X^+|)^2)^2}{|X^-|^2+\psi(|X^+|)^2}.$$ Then $F$ is a smooth function with $F(X) \geq 2\det X$ and $F(Du)=2 \det Du$. As before, this immediately implies that $u$ is a local minimiser of $\int_B F(Dv) \dd x$. It remains to verify that $H(X)=F(X)-|X|^2$ satisfies $|H'(X)| =o(|X|).$ To see this, note that $$H(X)=\frac{2\psi(|X^+|)^2(\psi(|X^+|)^2-3|X^-|^2)}{|X^-|^2+\psi(|X^+|)^2}.$$ Differentiating this, after a lengthy calculation one sees that $$|H'(X)| \lesssim \psi(|X^+|)(1+|\psi'(|X^+|)|).$$ Using the growth of $\psi$ and $\psi'$ this completes the proof. 
\end{proof}
We instead consider a different notion of first order approximation. We set $$H(z)=F(z)-|z|^p$$
and consider integrands $F$ which satisfy: 
\begin{equation}\label{eqn C^1 strong approximation dirichlet}
|H(z)-H(w)| \lesssim (1+|z|+|w|)^{p-2}|z-w|.
\end{equation}

\begin{lemma}\label{lemma C^1 approximation stated as a limsup in C^1 case}
If $F \in C^1(\R^{N \times n})$ then \eqref{eqn C^1 strong approximation dirichlet} is equivalent to
\begin{equation}\label{eqn C^1 approximation in the C^1 case}
\limsup_{|z| \rightarrow \infty}\frac{|F'(z)-pA_p(z)|}{|z|^{p-2}} <\infty.
\end{equation}
\end{lemma}
\begin{proof}
Simply observe if \eqref{eqn C^1 strong approximation dirichlet} holds then $$|H'(z) \cdot \xi| = \lim_{t \searrow 0} \frac{|H(z+t\xi)-H(z)|}{|t|} \lesssim (1+|z|)^{p-2}|\xi|$$ which means $|H'(z)| \lesssim (1+|z|)^{p-2}$ as required. Conversely, if this holds then $$|H(z)-H(w)| \lesssim |z-w| \int_0^1 (1+|w+t(z-w)|)^{p-2} \dd t$$
which is equivalent to $(1+|z|+|w|)^{p-2}|z-w|$, see the proof of Lemma \ref{lemma equivalences of auxiliary V functinos}.
\end{proof}
To obtain a $\BMO$ regularity result, we use the following fact from \cite[Remark 5.9]{Diening2012}.
\begin{lemma}\label{lem: BMO estimates for laplace equation general p}
Let $\Omega \subset \R^n$ be a bounded domain, $1<p<\infty$ and $u \in \WW^{1,p}_\mathrm{loc}(\Omega;\R^N)$.
If $$-\Div (A_p(Du+g))=0$$ where $g \in \BMO(\Omega)$, then $Du \in \BMO_\mathrm{loc}$.
\end{lemma}
Suppose now that $F \in C^1$ satisfies \eqref{eqn C^1 approximation in the C^1 case}. In this case, any minimiser $u$ satisfies the Euler--Lagrange equation $-\Div F'(Du) =0$. By defining $g(z)=A_p^{-1}(\frac{1}{p}F'(z))-z$ which is a continuous function, we can rewrite the Euler--Lagrange equation as $-\Div (A_p(Du+g(Du)))=0.$ As long as $g$ is bounded, Lemma \ref{lem: BMO estimates for laplace equation general p} then gives $Du \in \BMO_\mathrm{loc}$. We can still make this argument precise without assuming $F \in C^1$. 
 
\begin{theorem}\label{Theorem BMO regularity C^1 approximation for all p}
Let $1<p<\infty$, and suppose $u \in \WW^{1,p}(\Omega;\R^N)$ is a local minimiser of $\int_\Omega F(Dv) \dd x$ where $F$ is a continuous integrand satisfying \eqref{eqn C^1 strong approximation dirichlet}. Then $$Du \in \BMO_\mathrm{loc}(\Omega;\R^{N \times n}).$$
\end{theorem}
\begin{proof}
For any $\Omega' \Subset \Omega$ and $\varphi \in C_c^\infty(\Omega';\R^N)$, and any $t>0$ we have $$0 \leq \frac{1}{t} \int_{\Omega'}F(Du+tD\varphi)-F(Du) \dd x. $$
Since $F(z)=|z|^p+H(z)$, after applying the bound \eqref{eqn C^1 strong approximation dirichlet}, taking $t \searrow 0$, and considering $\pm \varphi$, it is not hard to see that $$\left|\int_{\Omega'}pA_p(Du) \cdot D\varphi \dd x \right| \leq C \int_{\Omega'}(1+2|Du|)^{p-2}|D\varphi| \dd x.$$
Thus, we can define a linear map $D\varphi \mapsto \int_{\Omega'}pA_p(Du)\cdot D\varphi$ which is bounded on the subspace of gradients in the weighted space $\LL^1(\Omega', \omega \dd x)$ where $\omega=(1+2|Du|)^{p-2}$. By Hahn--Banach, this map extends to the whole space, and so there is some $a \in \LL^\infty(\Omega', \omega \dd x)$ such that $$\int_{\Omega'}pA_p(Du) \cdot D\varphi \dd x =-\int_{\Omega'}a \cdot D\varphi \omega\dd x .$$ If we define $b=a \omega$, we see that $u$ satisfies the equation $$-\Div (pA_p(Du)+b)=0$$ in $\Omega'$, where $|b(x)| \lesssim (1+2|Du|)^{p-2}.$ If we now write $g=A_p^{-1}(A_p(Du)+\frac{1}{p}b)-Du$, we see that $$-\Div A_p(Du+g)=0.$$ By Lemma \ref{lem: BMO estimates for laplace equation general p}, we just need to show that $g \in \LL^\infty_\mathrm{loc}(\Omega)$. To see this, note that by taking inverses in Lemma \ref{lemma equivalences of auxiliary V functinos} we have
$$|A_p^{-1}(z)-A_p^{-1}(w)| \approx (|z|+|w|)^{\frac{2-p}{p-1}}|z-w|.$$
Applying this to $g$ gives
$|g| \approx (|A_p(Du)|+|b|)^{\frac{2-p}{p-1}}|b|.$
Since $|A_p(Du)|=|Du|^{p-1}$ and $|b| \lesssim (1+|Du|)^{p-2}$, it follows that $g$ is bounded, and hence we are done.
\end{proof}

By the same counterexample as in Theorem \ref{theorem counterexample non lipschitz extremely close}, it follows that Theorem \ref{Theorem BMO regularity C^1 approximation for all p} is sharp in the sense that Lipschitz regularity cannot be obtained in general. Indeed, $F(X)=|X|^2+(1-2|X^-|^2)_+$ satisfies \eqref{eqn C^1 strong approximation dirichlet} with $p=2$, simply because the function $t \mapsto (1-2t^2)_+$ is globally Lipschitz.

\section{Second-order asymptotic regularity}
In this section we consider integrands $F$ asymptotic to $|z|^p$ in a second order sense. For $F \in C^1(\R^{N \times n})$ we define
$$R_a(z)=\left[F(z)-F(a)-F'(a)\cdot (z-a)\right]-\left[|z|^p-|a|^p-pA_p(a)\cdot (z-a)\right].$$
In other words, $R_a(z)$ is the difference in first order Taylor expansions of $F$ and $|z|^p$ centred at $a$. We consider integrands $F \in C^1(\R^{N \times n})$ satisfying\footnote{We understand the fraction to equal $0$ whenever $z=a$.}
\begin{equation}\label{eqn approximate dirichlet in C^2 sense with taylor}
\lim_{|a| \rightarrow \infty} \sup_{z \in \R^{N \times n}}\frac{|R_a(z)|}{|V_p(z)-V_p(a)|^2}=0.
\end{equation}
Firstly to justify this assumption, we note that it forces $F$ to resemble $|z|^p$ in a $C^0$ sense. Indeed, if $F \in C^1(\R^{N \times n})$ satisfies \eqref{eqn approximate dirichlet in C^2 sense with taylor} then for every fixed $a \in \R^{N \times n}$, we have
$$\limsup_{|z|\rightarrow \infty}\frac{|F(z)-|z|^p|}{|z|^p}=\limsup_{|z| \rightarrow \infty}\frac{|R_a(z)|}{|V_p(z)-V_p(a)|^2} \leq \sup_{z \in \R^{N \times n}}\frac{|R_a(z)|}{|V_p(z)-V_p(a)|^2}$$ Taking the limit $a \rightarrow \infty$ gives $$\lim_{|z|\rightarrow \infty}\frac{|F(z)-|z|^p|}{|z|^p}=0.$$
We can equivalently characterise \eqref{eqn approximate dirichlet in C^2 sense with taylor} in the following way.
\begin{lemma}\label{lemma alternative form of C^2 approx}
Let $F \in C^1(\R^{N \times n})$. Then $F$ satisfies \eqref{eqn approximate dirichlet in C^2 sense with taylor} if and only if: for every $\ep>0$ there exists $C_\ep, M_\ep>0$ such that $$|R_a(z)| \leq \ep|V_p(z)-V_p(a)|^2+C_\ep \mbf{1}_{|a| \leq M_\ep}$$ for all $a, z \in \R^{N \times n}$. 
\end{lemma}
\begin{proof}
The statement trivially implies \eqref{eqn approximate dirichlet in C^2 sense with taylor}, so assume that \eqref{eqn approximate dirichlet in C^2 sense with taylor} holds. Given $\ep>0$, there exists $M_\ep>0$ such that $$\sup_{z \in \R^{N \times n}} \frac{|R_a(z)|}{|V_p(z)-V_p(a)|^2}<\ep$$ for all $|a| > M_\ep.$ So, we only need to consider $|a| \leq M_\ep$. Since $F(z)=|z|^p+o(|z|^p)$, and $F,F',A_p$ are all continuous, it follows that $R_a(z) = o(|z|^p)$ uniformly in $|a| \leq M_\ep$. So, there is some $R_\ep>0$ such that if $|z| > R_\ep$ and $|a| \leq M_\ep $ then $|R_a(z)| \leq \ep|V_p(z)-V_p(a)|^2$. Finally, if $|z| \leq R_\ep$ and $|a| \leq M_\ep$ we can simply estimate $|R_a(z)| \leq C_\ep$ by continuity. 
\end{proof}
If $F \in C^2(\R^{N \times n})$ then our condition is equivalent to the one considered in \cite{GiaquintaModica1986}.
\begin{lemma}
If $F \in C^{2}(\R^{N \times n})$ then \eqref{eqn approximate dirichlet in C^2 sense with taylor} is equivalent to
$$\lim_{|z| \rightarrow \infty}\frac{|F''(z)-D^2(|z|^p)|}{|z|^{p-2}}=0.$$
\end{lemma}
\begin{proof}
First assume that \eqref{eqn approximate dirichlet in C^2 sense with taylor} holds. Write $H(z)=F(z)-|z|^p$ as before, so that $R_a(z)=H(z)-H(a)-H'(a) \cdot (z-a).$ For $t>0$ and $|\xi|=1$ we have by Taylor's theorem (for large $|a|$)
$$R_a(a+t\xi)=t^2 \int_0^1 (1-s)H''(a+st\xi)\xi \cdot \xi \dd s = \frac{t^2}{2}H''(a) \xi \cdot\xi+o(t^2). $$
Moreover, by Lemma \ref{lemma equivalences of auxiliary V functinos} it follows that $|V_p(a+t\xi)-V_p(a)|^2 \approx t^2 |a|^{p-2}$ whenever $|t|<\frac{|a|}{2}$. Combining these facts we have for $|t|<\frac{|a|}{2}$ that
$$\frac{|H''(a)\xi \cdot \xi|}{|a|^{p-2}} \lesssim \frac{|R_a(a+t\xi)|+o(t^2)}{|V_p(a+t\xi)-V_p(a)|^2}.$$
Taking the limit $t \searrow 0$, for each fixed $a \neq 0$ we have
$$\frac{|H''(a) \xi \cdot \xi|}{|a|^{p-2}} \lesssim \lim_{t \searrow 0}\frac{|R_a(a+t\xi)|}{|V_p(a+t\xi)-V_p(a)|^2} \lesssim \sup_{z \in \R^{N \times n}} \frac{|R_a(z)|}{|V_p(z)-V_p(a)|^2}.$$
Taking $|a| \rightarrow \infty$, we deduce that $$\lim_{|a| \rightarrow \infty}\frac{|H''(a)|}{|a|^{p-2}}=0$$ as required. Conversely, suppose that this limit is zero. We let $\ep>0$ and take $|a|$ sufficiently large in terms of $\ep$. Then if $|z-a| \leq \frac{|a|}{2}$ it follows by Taylor's theorem
$$|R_a(z)| \lesssim |z-a|^2 \int_0^1(1-t)|H''(a+t(z-a)|) \dd t \lesssim \ep |z-a|^2 |a|^{p-2} \approx \ep |V_p(z)-V_p(a)|^2.$$ If instead $|z-a| > \frac{|a|}{2}$, so that $|z-a| \approx |a|+|z|$ and $|V_p(z)-V_p(a)|^2 \approx (|z|+|a|)^{p}$, we trivially bound 
\begin{align*}
    |R_a(z)| & \leq |H(z)|+|H(a)|+|H'(a)||z-a|\\
    & \leq \ep(|z|^p+|a|^p+|a|^{p-1}|z-a|)+C_\ep\\
    & \lesssim \ep |V_p(z)-V_p(a)|^2 +C_\ep
\end{align*}
which since $|V_p(z)-V_p(a)|^2 \gtrsim |a|^p$ in this case gives
$\frac{|R_a(z)|}{|V_p(z)-V_p(a)|^2} \lesssim \ep + \frac{C_\ep}{|a|^p}$.
Combining the two cases, we clearly have
$$\lim_{|a| \rightarrow \infty} \sup_{z \in \R^{N \times n}} \frac{|R_a(z)|}{|V_p(z)-V_p(a)|^2} \lesssim \ep$$ which gives condition \eqref{eqn approximate dirichlet in C^2 sense with taylor}.
\end{proof}

We now give a proof of the main result of this section. Our argument is based on  \cite{Diening2011, GiaquintaModica1986, ChipotEvans1986}, and the main difference here is that we do not assume that the integrand is $C^2$.

\begin{theorem}\label{theorem lipschitz regularity under approximate C^2 sense}
Suppose that $u \in \WW^{1,p}(\Omega;\R^N)$ is a local minimiser of $\int_\Omega F(Dv) \dd x$ where $F \in C^1(\R^{N \times n})$ satisfies \eqref{eqn approximate dirichlet in C^2 sense with taylor}. Then $u \in \WW^{1,\infty}_\mathrm{loc}(\Omega;\R^N)$ and $$\sup_{B_{R/2}}|Du|^p \ \leq C \left(1+\strokedint_{B_R}|Du|^p \dd x \right)$$
for all $B_R \Subset \Omega.$
\end{theorem}
\begin{proof}
Our aim will be to show the condition in Lemma \ref{lemma dichotomy boundedness} holds with $H=V_p(Du).$ To this end, we fix some ball $B \Subset \Omega$, and let $h \in u+\WW_0^{1,p}(B;\R^N)$ be $p$-harmonic. Using Lemma \ref{lemma equivalences of auxiliary V functinos} and the Euler--Lagrange equation for $h$, we have (as in the proof of Theorem \ref{theorem BMO in the case of extremely close})
$$\int_B|V_p(Du)-V_p(Dh)|^2 \dd x \approx \int_B |Du|^p -|Dh|^p \dd x.$$
By definition we have $$|Du|^p-|Dh|^p=R_a(Dh)-R_a(Du)+\left[F(Du)-F(Dh)\right]-F'(a)\cdot(Du-Dh)+pA_p(a)\cdot(Du-Dh).$$
It follows by minimality of $u$ and the divergence theorem that
$$\int_B|V_p(Du)-V_p(Dh)|^2 \dd x \lesssim \int_B R_a(Dh)-R_a(Du) \dd x.$$
Next, we set $\ep>0$ to be specified later. Recall that by Lemma \ref{lemma alternative form of C^2 approx}, there exists $C_\ep, M_\ep>0$ such that  $$|R_a(z)| \leq \ep|V_p(z)-V_p(a)|^2+C_\ep \mbf{1}_{|a| \leq M_\ep}$$ for all $a, z \in \R^{N \times n}$. Setting $V_p(a)=V_p(Du)_B$  this gives us 
\begin{align*}
\strokedint_B|V_p(Du)-V_p(Dh)|^2 \dd x &\leq \ep \strokedint_B|V_p(Dh)-V_p(Du)_B|^2+|V_p(Du)-V_p(Du)_B|^2\dd x \\
& \ + C_\ep \mbf{1}_{|a| \leq M_\ep}.
\end{align*}
Applying the triangle inequality, and ensuring $\ep$ is small enough that we can absorb the relevant term to the left hand side, this gives
$$\strokedint_B|V_p(Du)-V_p(Dh)|^2 \dd x \leq \ep \strokedint_B |V_p(Du)-V_p(Du)_B|^2 \dd x + C_\ep \mbf{1}_{|a| \leq M_\ep}.$$ Next, by running the exact same argument as in the proof of Theorem \ref{theorem BMO in the case of extremely close}, where we in particular use the decay estimates for $p$-harmonic functions, we have for all $\tau \in (0,\frac{1}{2})$
\begin{align*}
\strokedint_{\tau B}|V_p(Du)-V_p(Du)_{\tau B}|^2 \dd x &\leq C \tau^{2\alpha}\strokedint_{B}|V_p(Du)-V_p(Du)_B|^2\dd x\\
&+C\tau^{-n}\strokedint_{B}|V_p(Du)-V_p(Dh)|^2 \dd x
\end{align*}
where $\alpha \in (0,1)$ and $C>0$ depending only on $n,N,p$ and $F$. So, simply combining the above estimates we have
\begin{align*}
\strokedint_{\tau B}|V_p(Du)-V_p(Du)_{\tau B}|^2 \dd x &\leq (C \tau^{2\alpha}+C\tau^{-n}\ep)\strokedint_{B}|V_p(Du)-V_p(Du)_B|^2 \dd x\\
& + C\tau^{-n} C_\ep \mbf{1}_{|a| \leq M_\ep}.
\end{align*}
Now we may choose $\tau$ small enough so that $C\tau^{2\alpha}< \frac{1}{8}$, and then choose $\ep>0$ small enough so that $C\tau^{-n}\ep < \frac{1}{8}$. Then, our estimate reads
$$\strokedint_{\tau B}|V_p(Du)-V_p(Du)_{\tau B}|^2 \dd x \leq \frac{1}{4} \strokedint_B |V_p(Du)-V_p(Du)_B|^2 \dd x + C\mbf{1}_{|a| \leq M}$$
for some $C,M>0$. If $|a| >M$ we obtain the excess decay inequality in the dichotomy of Lemma \ref{lemma dichotomy boundedness}. Otherwise, $|a| <M$ and we have two cases. If $C < \frac{1}{4}\strokedint_B |V_p(Du)-V_p(Du)_B|^2 \dd x$ then once again we have the excess decay inequality, so we only need to consider the case $4C > \strokedint_B |V_p(Du)-V_p(Du)_B|^2 \dd x.$ But then recalling that $V_p(a)=V_p(Du)_B$ we have \begin{align*}
\strokedint_B|Du|^p \dd x &= \strokedint_B |V_p(Du)|^2 \dd x \\
& \leq \strokedint_B |V_p(Du)-V_p(Du)_B|^2 \dd x + |V_p(a)|^2 &\\
& \leq 4C+M^p := \tilde M
\end{align*}
which is the case in the dichotomy of Lemma \ref{lemma dichotomy boundedness}. In other words, we have shown that either $$\strokedint_{\tau B} |V_p(Du)-V_p(Du)_{\tau B}|^2 \dd x \leq \frac{1}{2}\strokedint_{B}|V_p(Du)-V_p(Du)_B|^2 \dd x $$ or $$\strokedint_B|Du|^p \dd x \leq \tilde M$$ for some constant $\tilde M$ depending only on $n,N,p$ and $F$. The conclusion of the theorem now follows by Lemma \ref{lemma dichotomy boundedness}. 
\end{proof} 
As a final remark, we note that there is no asymptotic condition that lets us deduce continuity of the gradient. 
\begin{prop}\label{prop counterexample to continuity of the gradient in C^2 approx}
There exists a smooth integrand $F:\R^{N \times n}\rightarrow \R$ which satisfies $F(z)=|z|^2$ for all $|z| \geq 2$, such that $\int_B F(Dv) \dd x$ has a local minimiser with discontinuous gradient. 
\end{prop}
\begin{proof}
For instance, we may take $F(z)=\eta(z)|z|^2$ where $\eta \in C^\infty(\R^{N \times n})$ is non-negative, vanishes at $z=\pm e_1 \otimes e_1$, and $\eta=1$ for $|z| \geq 2$. Then the map $u(x)=|x_1|e_1$ satisfies $F(Du)=0$. Since $F$ is non-negative, this means $u$ is a local minimiser. However clearly $Du$ is not continuous.
\end{proof}

\section{Appendix: BMO functions}
We collect some basic facts about functions of bounded mean oscillation, see \cite{JohnNirenberg1961, Grafakos2014}.
We define $$\BMO(\R^n)=\left\{u \in \LL^1_\mathrm{loc}(\R^n): [u]_{\BMO}=\sup_{Q \subset \R^n} \strokedint_Q |u-(u)_Q| \dd x <\infty  \right\}$$
The following theorem due to John--Nirenberg is well known. 
\begin{theorem}
Let $u \in \LL^1_\mathrm{loc}(\R^n).$ Then the following are equivalent:
\begin{enumerate}[(i)]
    \item $u \in \BMO (\R^n)$.
    \item There exists $c_1, c_2 >0$ such that $$ \sup_{Q \subset \R^n} \frac{|\left\{x \in Q: |u(x)-(u)_Q|>t \right\}|}{|Q|} \leq c_2 e^{-c_1 t}.$$
    \item There exists $\alpha >0$ such that $$\sup_{Q \subset \R^n}\strokedint_Q e^{\alpha|u-(u)_Q|} \dd x <\infty.$$
\end{enumerate}
Consequently, $\BMO(\R^n) \hookrightarrow \LL^p_\mathrm{loc}(\R^n)$ for all $1 \leq p <\infty$. Moreover, $$[u]_\BMO \approx \sup_{Q \subset \R^n} \left(\strokedint_Q|u-(u)_Q|^p \dd x \right)^{1/p}.$$
\end{theorem}
One can equivalently define $\BMO$ using balls instead of cubes. The following simple but useful lemma gives a condition to show a map is not in $\BMO$.

\begin{lemma}\label{lemma function not in BMO}
Suppose that on some fixed ball $B_r$, we have that $$|\left\{x \in B_r: |u(x)|>t \right\}| \geq c_1 e^{-c_2 {t^{1/\alpha}}}$$ for some $c_1, c_2 >0$ and $\alpha >1$, and for all $t$ sufficiently large. Then $u \notin \BMO$.
\end{lemma}
\begin{proof}
Simply observe that if $t > 2 |(u)_{B_r}|$ then $$\left\{x \in B_r: |u(x)| > t\right\} \subset \left\{x \in B_r: |u(x)-(u)_{B_r}| > t/2 \right\}.$$ Therefore, if the claimed inequality holds, by the John--Nirenberg lemma we must have that $c_1 e^{-c_2 {t^{1/\alpha}}} \leq c_3 e^{-c_4 t}$ for all large $t$, which is impossible.
\end{proof}
Roughly speaking, $\BMO$ maps are allowed to have logarithmic singularities, but singularities such as $|\log |x||^\alpha$ where $\alpha>1$ are not allowed.

\section{AI Declaration}
The author benefitted from the use of ChatGPT-6 Astra in the development of this paper. The idea to employ the direct method to prove regularity is classical, and the precise way we have implemented it here is a combination of the author's and ChatGPT's ideas. Similarly, the idea to use conformal splitting to produce counterexamples is well known in the literature, and the author combined this with suggestions from ChatGPT in the present examples. All text and equations were written by hand solely by the author.

\end{document}